\documentclass[11pt,twoside]{article}

\usepackage{a4wide}
\usepackage{amsfonts}
\usepackage{amssymb}
\usepackage{amsmath}
\usepackage{graphicx}

\newcommand{\ignore}[1]{}

\def\@begintheorem#1#2{\par\bgroup{\sc #1\ #2. }\it\ignorespaces}
\def\@opargbegintheorem#1#2#3{\par\bgroup{\sc #1\ #2\ (#3). } \it\ignorespaces}
\def\@endtheorem{\egroup}
\newtheorem{theorem}{Theorem}[section]
\newtheorem{corollary}[theorem]{Corollary}
\newtheorem{lemma}[theorem]{Lemma}

\newtheorem{example}[theorem]{Example}
\newtheorem{proposition}[theorem]{Proposition}
\newtheorem{definition}[theorem]{Definition}
\newcommand{\bt}[1]{\begin{theorem}\label{#1}}
\newcommand{\bc}[1]{\begin{corollary}\label{#1}}
\newcommand{\bl}[1]{\begin{lemma}\label{#1}}
\newcommand{\be}[1]{\begin{example}\label{#1}}
\newcommand{\bp}[1]{\begin{proposition}\label{#1}}
\newcommand{\ba}[1]{\begin{algorithm}\rm\label{#1}}
\newcommand{\bd}[1]{\begin{definition}\rm\label{#1}}{\normalsize }
\newcommand{\bpr}{\noindent {\em Proof. }}
\newcommand{\et}{\end{theorem}}
\newcommand{\ec}{\end{corollary}}
\newcommand{\el}{\end{lemma}}
\newcommand{\ee}{\end{example}}
\newcommand{\ep}{\end{proposition}}
\newcommand{\ed}{\end{definition}}
\newcommand{\epr}{{\ \vbox{\hrule\hbox{%
\vrule height1.3ex\hskip0.8ex\vrule}\hrule}}\\\par}

\def\R{\mathbb{R}}
\def\Z{\mathbb{Z}}
\def \conv {{\rm conv}}

\begin{document}

\title{\bf On Degree Sequence Optimization}

\author{
Shmuel Onn
\thanks{\small Technion - Israel Institute of Technology.
Email: onn@ie.technion.ac.il}
}
\date{}

\maketitle

\begin{abstract}
We consider the problem of finding a subgraph of a given graph which maximizes a given function evaluated at its degree sequence. While the problem is intractable already for convex functions, we show that it can be solved in polynomial time for convex multi-criteria objectives.
We next consider the problem with separable objectives, which
is NP-hard already when all vertex functions are the square. We consider a colored extension of the separable problem, which includes the notorious exact matching problem as a special case, and show that it can be solved in polynomial time on graphs of bounded tree-depth for any vertex functions. We mention
some of the many remaining open problems.

\end{abstract}

\section{Introduction}

The {\em degree sequence} of a simple graph $H=(V,E)$ with $V=[n]:=\{1,\dots,n\}$
is the vector $d(H)=(d_1(H),\dots,d_n(H))$, where $d_i(H):=|\{e\in E:i\in e\}|$
is the degree of vertex $i$ for all $i\in[n]$. If $G=([n],F)\subseteq H$ with $F\subseteq E$
is a subgraph of $H$ we also write $d(F):=d(G)$ for the degree sequence of $G$.
Degree sequences have been studied by many authors, starting from their characterization
by Erd\H{o}s and Gallai \cite{EG}, see e.g. \cite{EKM} and the references therein.

\vskip.2cm
In this article we are interested in the following problem and some of its variants.

\vskip.2cm\noindent{\bf Degree Sequence Optimization.}
Given a graph $H=([n],E)$, integer $0\leq m\leq|E|$, and function $f:\Z^n\rightarrow\Z$,
find a subgraph $G=([n],F)\subseteq H$ with $m$ edges maximizing $f(d(G))$.

\vskip.2cm
We will also consider the
{\em unprescribed} variant of the problem, where the number of edges $m$ is not specified, and the optimization is over any subgraph
$G\subseteq H$ with vertex set $[n]$.

\vskip.2cm
Our subgraphs always consist of the entire original
vertex set and a subset of the edges.
Scaling up rational values if needed we assume the function takes on integer values.
Also, it is enough that the function
$f$ is defined only on its domain $\{0,1,\dots,d_1(H)\}\times\{0,1,\dots,d_n(H)\}$,
and properties of the function such as convexity are with respect to this domain only.
We assume that the function is presented either by an oracle that, queried on $z\in\Z^n$,
returns $f(z)$, or by some compact presentation as will be clear from the context.

As an example, consider the case where the function is {\em linear} and given by the inner
product $f(x)=wx$ for some $w\in\Z^n$. Then for every subset of edges $F\subseteq E$ we have
that $wd(F)=\sum\{wd(e):e\in F\}$ and so the problem can be easily solved by sorting the edges
of $E$ by the value $wd(e)$ and taking $F$ to consist of the $m$ edges with maximum values.

However, the problem is generally much harder, even for convex functions, and when such
functions are presented by oracles, any algorithm for the problem may need to make
exponentially many queries and hence the problem is intractable, as explained in Section 4.

In contrast to this intractability, and as a natural extension of the linear case,
we show that for the following convex multi-criteria objective the problem is
polynomial time solvable. Let $w_1,\dots,w_r\in\Z^n$ be given vectors. Each $w_i$ is interpreted as a
linear criterion under which the value of subgraph $G\subseteq H$ is the inner product $w_id(G)$.
These criteria are balanced by a convex function $f:\Z^r\rightarrow\Z$. The case of a
single criterion $r=1$ and $f$ the identity on $\Z$ is the linear problem discussed above.
Our first result is the efficient solution of this problem. Here and elsewhere an algorithm
involving oracle presented objects is {\em polynomial time} if the number of arithmetic
operations and oracle queries it makes are polynomial in the data.
\bt{multicriteria_m}
Fix any $r$. Given a graph $H=([n],E)$, $m\in\Z_+$, $w_1,\dots,w_r\in\Z^n$, and oracle presented
convex $f:\Z^r\rightarrow\Z$, we can solve in polynomial time the multi-criteria problem
$$\max\{f\left(w_1 d(G),\dots,w_r d(G)\right)\ :\ G=([n],F)\subseteq H,\ \ |F|=m\}\ .$$
\et

Next we consider {\em separable} functions, that is, of the form
$f(x)=\sum_{i=1}^n f_i(x_i)$ with each $f_i:\Z\rightarrow\Z$ univariate.
In this case several results are known in the literature. First, even in this case, the problem
is generally NP-hard, see
Section 4. On the positive side,
the problem is polynomial time solvable in the following situations: over the complete graph $H=K_n$ when $f_i(z)=z^2$ for all $i$ \cite{PPS} and more generally when all functions are the same $f_1=\cdots=f_n$ and arbitrary
\cite{DLMO}; for any graph $H$ when all $f_i$ are concave for the prescribed problem \cite{AS} and independently for the unprescribed variant \cite{DO}; and for every fixed $r$ and any $H$ when $r$ of the functions are arbitrary and the rest are either all nondecreasing or all nonincreasing \cite{DO}.

A further special case of the separable problem is the {\em general factor problem} \cite{Cor},
which is to decide, given a graph $H$ and subsets $B_i\subset\Z$ for $i\in[n]$,
if there is a $G\subseteq H$ with $d_i(G)\in B_i$ for all $i$, and find one if yes.
Indeed, for each $i$ define a function $f_i$ by $f_i(x):=0$ if $x\in B_i$ and $f_i(x):=-1$
if $x\notin B_i$. Then the optimal value of the degree sequence problem is zero if and only if
there is a factor, in which case any optimal graph $G$ is one. An even more special case
is the well studied {\em $(l,u)$-factor problem} introduced by Lov\'asz \cite{Lov},
where each $B_i=\{l_i,\dots,u_i\}$ is an interval. This reduces to the degree sequence
problem even with concave functions, with $f_i(z):=z-l_i$ if $z\leq l_i$, $f_i(z):=0$ if
$l_i\leq z\leq u_i$, and $f_i(z):=u_i-z$ if $u_i\leq z$. In particular, the
{\em perfect matching problem} is that with $B_i=\{1\}$ for all $i$.

We consider the following extension of the separable problem. We are now given,
with the graph $H=([n],E)$, a $p$-partition $E=\biguplus_{k=1}^pE_k$,
and integers $0\leq m_k\leq|E_k|$ for $k=1,\dots,p$. We search for a subgraph
$G=([n],F)\subseteq H$ maximizing the objective $\sum_{i=1}^n f_i(d_i(G))$ with the
requirement that $|F\cap E_k|=m_k$ for all $k$. We refer to the $E_k$ as {\em colors} and call the problem the {\em colored degree optimization problem}. Standard separable degree optimization is the special case of $p=1$.
Another special case is the notorious {\em exact matching problem}, where we are given
a $p$-partition $E=\biguplus_{k=1}^pE_k$ of the edges of the complete bipartite graph
$K_{r,r}$ and integers $m_1,\dots,m_p$, and we need to decide if there is a perfect matching
$M\subset E$ with $|F\cap E_k|=m_k$ for all $k$. It is a special case of our problem
with $H=K_{r,r}$ and $f_i(z)=-(z-1)^2$ for every vertex $i$, where the optimal value
is zero if and only if there is an exact matching.

The exact matching problem has randomized algorithms \cite{MVV,Onn} but its deterministic
complexity is open already for $p=3$. Moreover, standard separable degree optimization,
with $p=1$, is generally NP-hard, see Section 4. In contrast, we show that for graphs $H$ of
bounded tree-depth (see Section 3), the colored problem can be solved in polynomial time.

\bt{tree-depth}
Fix any $p,d$. Given a graph $H=([n],E)$ of tree-depth at most $d$,
functions $f_1,\dots,f_n:\Z\rightarrow\Z$,
$p$-partition $E=\biguplus_{k=1}^pE_k$, and integers $0\leq m_k\leq|E_k|$ for $k=1,\dots,p$,
we can solve in polynomial time the following colored degree sequence optimization problem
$$\max\left\{\sum_{i=1}^n f_i(d_i(G))\ :
\ G=([n],F)\subseteq H,\ \ |F\cap E_k|=m_k,\ \ k=1,\dots p\right\}\ .$$
\et

We conclude with some open problems. The complexity of degree sequence optimization is still wide open, even in the separable case. It would be interesting to determine it for various classes
of graphs $H$ and various classes of functions $f$. Particularly intriguing is the special case of the
complete graph $H=K_n$, where the optimization is over any graph $G$ on $[n]$, for which the separable problem
might be solvable in polynomial time for any functions $f_i$. Also, while the separable problem over arbitrary $H$ was shown to be NP-hard in \cite{AS} already when $f_i(z)=z^2$ for all $i$,
the complexity of the unprescribed variant, with no restriction on the number $m$ of edges,
is open and might be solvable in polynomial time for any graph $H$ and any convex functions $f_i$.
In particular, what is the complexity of the following unprescribed separable problem over the complete graph $H=K_n$
with arbitrary convex functions,
$$\max\left\{\sum_{i=1}^n f_i(d_i(G))\ :\ G\ \mbox{any graph on}\ [n],
\ \ f_1,\dots,f_n\ \mbox{any convex functions}\right\}\ ?$$

\vskip.2cm
In Section 2 we prove Theorem \ref{multicriteria_m}. In Section 3 we prove Theorem \ref{tree-depth}.
We conclude in Section 4 with some limitations on the solvability of the problem and its variants.

\section{Convex multi-criteria objectives}

We need some terminology. Let $S\subset\Z^n$ be a finite set and let $P=\conv(S)$.
Consider any edge ($1$-dimensional face) $C=[u,v]$ of $P$.
A {\em direction} of $C$ is any nonzero multiple of $v-u$.
A {\em linear optimization oracle} for $S$ is one that, queried on $w\in\Z^n$,
returns an element $x^*\in S$ attaining $\max\{wx:x\in S\}$.
We make use of the following result \cite[Theorem 2.16]{Onn}.
\bp{edge_directions}
Fix any $r$. Given any $S\subset\Z^n$ presented by a linear optimization oracle,
$w_1,\dots,w_r\in\Z^n$, oracle presented convex function $f:\Z^r\rightarrow\Z$, and a set
$D$ of directions of all edges of $P=\conv(S)$, we can solve in polynomial time
the multi-criteria problem
$$\max\{f\left(w_1 x,\dots,w_r x\right)\ :\ x\in S\}\ .$$
\ep

We can now prove our first theorem.

\vskip.2cm\noindent{\bf Theorem \ref{multicriteria_m}}
{\em Fix any $r$. Given a graph $H=([n],E)$, $m\in\Z_+$, $w_1,\dots,w_r\in\Z^n$, and oracle
presented convex $f:\Z^r\rightarrow\Z$, we can solve in polynomial time the multi-criteria problem
$$\max\{f\left(w_1 d(G),\dots,w_r d(G)\right)\ :\ G=([n],F)\subseteq H,\ \ |F|=m\}\ .$$
}

\bpr
Recall that for any subset $F\subseteq E$ we write $d(F):=d(G)$ for the degree sequence of the
subgraph $G=([n],F)\subseteq H$ with set of edges $F$. In particular, if $e=\{i,j\}\in E$ then
$d(e)={\bf 1}_i+{\bf 1}_j$ is the sum of the unit vectors corresponding to the vertices of $e$.
Define
$$S_H^m\ :=\ \{d(F)\ :\ F\subseteq E,\ \ |F|=m\}\ \subset\ \Z^n
\ ,\quad\quad\quad P_H^m\ :=\ \conv(S)\ \subset\ \R^n\ .$$

Assume $m<|E|$ else the problem is trivial. First we construct a set $D$ of directions
of every edge of $P_H^m$ of size polynomial in $n$.
Consider any edge $C$ of $P_H^m$ and let $d(F_1),d(F_2)$ be its vertices
so $C=[d(F_1),d(F_2)]$ with $F_1,F_2\subset E$ and $|F_1|=|F_2|=m$. Let $w\in\Z^n$ be such that
the inner product $wx$ is maximized over $P_H^m$ precisely at $C$. Since $|F_1|=m=|F_2|$
there exist $e_1\in F_1\setminus F_2$ and $e_2\in F_2\setminus F_1$.
If $w d(e_1)-w d(e_2)\geq 0$ then $w d(F_2\setminus e_2\cup\{e_1\})\geq w d(F_2)$ so
$d(F_2\setminus e_2\cup\{e_1\})\in C$ and hence
$d(F_2)-d(F_2\setminus e_2\cup\{e_1\})=d(e_2)-d(e_1)$ is a direction of the edge $C$.
Similarly if $w d(e_1)-w d(e_2)\leq 0$ then $d(F_1)-d(F_1\setminus e_1\cup\{e_2\})=d(e_1)-d(e_2)$
is a direction of $C$. Since directions are defined up to nonzero scalar multiplication, it follows
that the set $D$ consisting of one representative of $\pm(d(e)-d(f))$
for each pair of distinct $e,f\in E$, is a set of ${|E|\choose 2}=O(n^4)$
vectors containing a direction of every edge of $P_H^m$.

Next, we construct a linear optimization oracle for $S_H^m$. Given $w\in\Z^n$, we need to solve
$$\max\{wx:x\in S_H^m\}\ =\ \max\{wd(F)\ :\ F\subseteq E,\ \ |F|=m\}\ .$$
Now we observe that for every $F\subseteq E$ we have $wd(F)=\sum\{wd(e):e\in F\}$.
Thus, we let $F^*$ consist of the $m$ edges $e\in E$ with largest values $wd(e)$
and return $x^*:=d(F^*)$. Using this oracle and the set $D$ in
Proposition \ref{edge_directions} we obtain the claim of the theorem.
\epr
To solve the unprescribed variant of the problem where the number of edges is not specified
we can simply solve the above problem for $m=0,1,\dots,|E|$ and pick the best solution.
However, for the this variant there is a faster shortcut which we proceed to describe. Let now
$$S_H\ :=\ \{d(F)\ :\ F\subseteq E\}\ \subset\ \Z^n
\ ,\quad\quad\quad P_H\ :=\ \conv(S)\ \subset\ \R^n\ .$$
Consider any edge $C=[d(F_1),d(F_2)]$ of $P_H$ and let $w\in\Z^n$ be such that $wx$ is maximized
over $P_H$ precisely at $C$. Suppose first there exist $e\in F_1\setminus F_2$.
If $wd(e)\geq 0$ then we have $w d(F_2\cup\{e\})\geq w d(F_2)$ so $d(F_2\cup\{e\})\in C$ and
hence $d(F_2\cup\{e\})-d(F_2)=d(e)$ is a direction of the edge $C$. Likewise if
$w d(e)\leq 0$ then $d(F_1)-d(F_1\setminus e)=d(e)$ is a direction of $C$. A similar argument
works if there exist $e\in F_2\setminus F_1$. It follows that the set
$D:=\{d(e):e\in E\}$ is a set of $|E|=O(n^2)$ directions of every edge of $P_H$.
Next, given $w\in\Z^n$, we need to solve
$$\max\{wx:x\in S_H\}\ =\ \max\{wd(F)\ :\ F\subseteq E\}\ .$$
Since $wd(F)=\sum\{wd(e):e\in F\}$ we set $F^*:=\{e\in E:wd(e)>0\}$ and return $x^*:=d(F^*)$.

So now the set $D$ has size $O(n^2)$ instead of $O(n^4)$ in the prescribed variant,
and linear optimization can be done by checking signs in the set $\{wd(e):e\in E\}$
rather than sorting it in the prescribed variant.
So applying Proposition \ref{edge_directions} again we obtain a faster solution.

\section{Colored bounded tree-depth graphs}

We need some more terminology.
The tree-depth of a graph $G=(V,E)$ is defined as follows. The {\em height} of a rooted
tree is the maximum number of vertices on a path from the root to a leaf. A rooted tree
on $V$ is {\em valid} for $G$ if for each edge $\{i,j\}\in E$ one of $i,j$ lies on the
path from the root to the other of $i,j$. The {\em tree-depth} $td(G)$ of $G$ is the smallest height
of a rooted tree which is valid for $G$. For instance, if $G=([2m],E)$ is a perfect matching
with $E=\{\{i,m+i\}:i\in[m]\}$ then its tree-depth is $3$ where a
tree validating it rooted at $1$ has edge set $E\uplus\{\{1,i\}:i=2,\dots,m\}$.
Next, the {\em graph} of an $m\times n$ matrix $A$ is the graph $G(A)$ on $[n]$ where ${j,k}$ is
an edge if and only if there is an $i\in[m]$ such that $A_{i,j}A_{i,k}\neq 0$.
The {\em tree-depth} of $A$ is the tree-depth $td(A):=td(G(A))$ of its graph.
We use a recent result of \cite{EHKKLO,KLO} on integer programs in variable dimension $n$
(as opposed to the classical result in fixed dimension \cite{Len}). It asserts that
integer programming is solvable in fixed-parameter tractable time \cite{DF} when parameterized
by the {\em numeric measure} $a:=\|A\|_\infty:=\max_{i,j}|A_{i,j}|$ and {\em sparsity measure} $d:=td(A^T)$ of
the matrix defining the program. In particular, for any fixed $d$ it is solvable in
polynomial time even when $a$ is a variable part of the input and given in unary.

\bp{IP}
Consider the following integer programming problem in variable dimension $n$,
where $A\in\Z^{m\times n}$, $c,l,u\in\Z^n$, $b\in\Z^m$, parameterized by
$a:=\|A\|_\infty$ and $d:=td(A^T)$,
$$\max\{cx\ :\ Ax=b,\ l\leq x\leq u,\ x\in\Z^n\}\ .$$
It can be solved in time $(2a+1)^{O(d2^d)}n^2L$ where $L$ is the total bit size
of the data $A,b,c,l,u$.
\ep

We can now prove our second theorem.

\vskip.2cm\noindent{\bf Theorem \ref{tree-depth}}
{\em Fix any $p,d$. Given a graph $H=([n],E)$ of tree-depth at most $d$,
functions $f_1,\dots,f_n:\Z\rightarrow\Z$,
$p$-partition $E=\biguplus_{k=1}^pE_k$, and integers $0\leq m_k\leq|E_k|$ for $k=1,\dots,p$,
we can solve in polynomial time the following colored degree sequence optimization problem
$$\max\left\{\sum_{i=1}^n f_i(d_i(G))\ :
\ G=([n],F)\subseteq H,\ \ |F\cap E_k|=m_k,\ \ k=1,\dots p\right\}\ .$$
}
\bpr
We construct the following integer program, in $n+3|E|$ binary variables,
variable $x_e$ for each edge $e\in E$ and variable $y_i^j$ for each vertex
$i\in[n]$ and each $j\in\{0,\dots,d_i(H)\}$,
where, as usual, $\delta_i(H):=\{e\in E:i\in e\}$ is the set of edges in $H$ containing vertex $i$,
$$\max\sum_{i=1}^n\sum_{j=0}^{d_i(H)}f_i(j)\cdot y_i^j$$
\begin{equation}\label{IP1}
\sum_{e\in\delta_i(H)}x_e-\sum_{j=0}^{d_i(H)}j\cdot y_i^j\ =\ 0\ ,\quad\quad i\in[n]\ ,
\end{equation}
\begin{equation}\label{IP2}
\sum_{j=0}^{d_i(H)}y_i^j\ =\ 1\ ,\quad\quad i\in[n]\ ,
\end{equation}
\begin{equation}\label{IP3}
\sum_{e\in E_k}x_e\ =\ m_k\ ,\quad\quad k\in[p]\ ,
\end{equation}
$$x_e\in\{0,1\}\ ,\quad e\in E\ ,\quad\quad y_i^j\in\{0,1\}
\ ,\quad i\in[n]\ ,\quad j\in\{0,\dots,d_i(H)\}\ .$$
Suppose $(x,y)$ is a feasible solution of this program and let $G=([n],F)$ be the subgraph with
$F:=\{e\in E:x_e=1\}$. Consider any $i\in[n]$. Constraint \eqref{IP2} forces $y_i^j=1$
for exactly one $j\in\{0,\dots,d_i(H)\}$. Then constraint \eqref{IP1} forces this $j$ to be
$j=d_i(G)$. So the objective value of $(x,y)$ in the program is $\sum_{i=1}^n f_i(d_i(G))$
which is the objective value of $G$ in the degree optimization problem. Finally,
constraint \eqref{IP3} forces $|F\cap E_k|=m_k$ for all $k$.

It is easy to see that also, conversely, if $G=([n],F)$ is a feasible solution of the
degree optimization problem then $(x,y)$ defined by $x_e=1$ if $e\in F$ and $x_e=0$ otherwise,
and $y_i^j=1$ if $j=d_i(G)$ and $y_i^j=0$ otherwise,
is a feasible solution to the program with the same objective value.
So the degree optimization problem reduces to solving the integer program.

Consider the matrix $A$ expressing equations \eqref{IP1}--\eqref{IP3} and its transpose $A^T$.
The columns of $A$ are indexed by the variables. Let us index the equations and the rows of
$A$ by $a_1,\dots,a_n$ corresponding to equations \eqref{IP1}, $b_1,\dots,b_n$ corresponding
to equations \eqref{IP2}, and $c_1,\dots,c_p$ corresponding to equations \eqref{IP3}.
Let $T$ be a rooted tree on $\{a_1,\dots,a_n\}\cong[n]$ validating that $td(H)\leq d$
and let $a_r$ be its root. We now use $T$ to obtain a rooted tree $T'$
with vertices $a_1,\dots,a_n,b_1,\dots,b_n,c_1,\dots,c_p$, rooted at $c_1$, consisting of
the edges of $T$, the edges $\{a_i,b_i\}$ for $i\in[n]$, the edges $\{c_i,c_{i+1}\}$
for $1\leq i<p$, and the edge $\{c_p,a_r\}$.

We now show that $T'$ is valid for $G(A^T)$.
For this we need to show that if two equations share a variable then
one lies on the path in $T'$ from its root to the other. Consider any $a_s,a_t$.
They share a variable $x_e$ if only if $e=\{s,t\}\in E$ and then, since $T$ is valid for $H$,
one of them lies on the path in $T$ from its root $a_r$ to the other, and hence also on the
path in $T'$ from its root $c_1$ to the other. Next, any distinct $b_s,b_t$ do not share a variable.
Next, consider any $a_s,b_t$. They share the variables $y_s^j=y_t^j$ if and only if $s=t$ in which case
$a_s$ lies on the path in $T'$ from its root to $b_t$ consisting of the path from $c_1$ to $a_s$ via
$c_p$ and $a_r$ and the edge $\{a_s,b_t\}$. Next, any distinct $c_s,c_t$ do not share a variable and any $c_s,b_t$ do not share a variable. Finally, for any $c_s,a_t$ we have that $c_s$ lies on the path from $c_1$ to $a_t$ . So $T'$ is valid for $G(A^T)$.

Now, the height of $T'$ is at most $d':=p+d+1$. Also, the matrix satisfies $a=\|A\|_\infty\leq n-1$
attained by equations \eqref{IP1}. The total bit size $L$ of $A,b,c,l,u$ is clearly polynomial in
the data since the entries of $A$ are bounded by $a\leq n-1$, the entries of $b$ are bounded
by the $m_k$, the entries of $c$ come from the function values $f_i(j)$, the entries of $l$ are
all $0$ and the entries of $u$ are all $1$. By Proposition \ref{IP}, and since $p,d$ are fixed,
the integer program and hence the colored degree problem are solvable
in polynomial time $(2a+1)^{O(d'2^{d'})}n^2L=n^{O(d'2^{d'})}L$.
\epr

\section{Some limitations}
\label{hem}

Here we discuss some situations where the degree sequence optimization problem is hard.

\subsubsection*{Oracle presented convex functions}

Given a graph $H=([n],E)$ consider the {\em degree sequence polytope} of $H$ defined as
$$P_H\ :=\ \conv\{d(G)\ :\ G\subseteq H\}\ \subset\ \R^n\ .$$
When maximizing a convex $f$ there will be an optimal graph $G$ with $d(G)$
a vertex of $P_H$. If $f$ is presented by an oracle then we claim that any algorithm
for the problem may need to make exponentially many queries and hence the problem
is intractable. To see this, consider a perfect matching graph $H=([2m],E)$ with $n=2m$
vertices and $E=\{\{i,m+i\}:i\in[m]\}$. Then $P_H$ is isomorphic to the
unit cube $[0,1]^m\subset\R^m$ via the map $[0,1]^m\rightarrow P_H:x\mapsto(x,x)$.
Any integer values on the vertices of $P_H$ extend to a convex function
$f:\Z^n\rightarrow\Z$ and hence any algorithm for the problem
must query the oracle on each of the $2^m$ vertices of $P_H$.

\subsubsection*{Concave-convex separable functions}

We now show the hardness of the unprescribed variant of the problem with separable functions.
The NP-complete {\em cubic subgraph problem} is to decide if a given graph
$H=([n],E)$ has a subgraph $G$ where each vertex has degree $0$ or $3$.
Defining $f_i(0)=f_i(3)=0$ and $f_i(z)=-1$ for $z\neq0,3$ for $i=1,\dots,n$,
the optimal objective value is zero if and only if $H$ has a cubic subgraph,
and so the corresponding degree optimization problem is NP-hard.

The problem remains NP-hard with $H=(I,J,E)$ bipartite and $f_i$ concave for
$i\in I$ and convex for $i\in J$. Indeed, the general factor problem is NP-complete
for bipartite $H=(I,J,E)$ with maximum degree $3$ and $B_i=\{1\}$ for $i\in I$
and $B_i=\{0,3\}$ for $i\in J$, see \cite{Cor}. Define
$$
f_i(z)\ :=\ -(z-1)^2
\quad i\in I,
\quad\quad\quad\quad
f_i(z)\ :=\ z(z-3)
\quad i\in J\quad.
$$
Then the optimal value of the degree sequence problem is zero if and only if there is a factor.

It remains hard moreover for bipartite graphs with a single concave function and all others convex. Recall that the prescribed problem on $H$ with specified number $m$ of edges and all functions $f_i(z)=z^2$ is NP-hard \cite{AS}. Define a bipartite graph
$L$ by subdividing each edge $\{i,j\}$ of $H$ and denoting the new
vertex by $\{i,j\}$, and adding a new vertex $s$ connected to all
$\{i,j\}$ vertices. Define $f_i(z):=z^2$ for all original vertices, and, for a sufficiently large positive integer $a$, let $f_s(z):=-a(z-m)^2$ and $f_{i,j}(z):=az(z-3)$ for all
new vertices $\{i,j\}$. Then in
any optimal subgraph $G\subseteq L$
for the unprescribed problem on $L$,
$m$ of the vertices
$\{i,j\}$ have degree $3$ and the rest have degree $0$, and the subgraph of $H$ with edge set $\{\{i,j\}:d_{i,j}(G)=3\}$ is optimal for the prescribed problem
on $H$, reducing the latter to the former.

\subsubsection*{Weighted degree optimization}

Another extension is the following problem. With the graph
$H=([n],E)$ we are given an edge weighting $w:E\rightarrow\Z$.
The problem is to find a subgraph $G=([n],F)\subseteq H$ maximizing
$$\sum_{i=1}^n f_i\left(\sum\{w(e):e\in\delta_i(G)\}\right)\ .$$
where $\delta_i(G)=\{e\in F:i\in e\}$ is the set of edges in $G$ containing $i$.
So the standard separable degree sequence optimization problem is the
special case with $w(e)=1$ for all $e\in E$.

However, this is NP-hard already for $H=(\{v_1,v_2\},[m],E)\simeq K_{2,m}$ with concave
$f_{v_1}(z):=-(z-a)^2$ and $f_{v_2},f_1,\dots,f_m$ all zero. Recall that the NP-complete
{\em partition problem} is to decide, given positive integers $a_1,\dots,a_m$, if there
is a $J\subset[m]$ with $\sum_{j\in J}a_j=a:={1\over2}\sum_{j=1}^ ma_j$.
Given such integers, define the weight function $w:E\rightarrow\Z$ by $w(\{v_i,j\}):=a_j$
for each $i=1,2$ and $j\in[m]$, and define the functions as above.
Then clearly there is a partition if and only if the optimal objective value of the degree sequence problem is zero, showing its hardness.

\section*{Acknowledgments}
S. Onn was supported by a grant from the Israel Science Foundation and the Dresner chair.

\end{document}